\input amstex
\raggedbottom
\documentstyle{amsppt}
\magnification=1200
%\hcorrection{-0.3cm}
%\vcorrection{-0.3cm}
\pageheight{9.0 true in} \pagewidth{6.5 true in}
%\nologo
\pageno=1 \NoRunningHeads
\topmatter
\title
 $q$-Volkenborn Integration and Its Applications \endtitle
\author
    TAEKYUN KIM \endauthor
 \affil
  Department of Mathematics Educations,
   Kongju National University,
   Kongju 314-701, S. Korea\\
    $\text{e-mail: tkim$\@$kongju.ac.kr}$\\
  \endaffil
 \dedicatory
 \enddedicatory
 \abstract
  The main purpose of this paper is to present a systemic study of
  some families of multiple $q$-Euler numbers and polynomials.
  In particular, by using the $q$-Volkenborn  integration on $\Bbb
  Z_p$, we construct $p$-adic $q$-Euler numbers and polynomials of
  higher order. We also define new generating functions of
  multiple $q$-Euler numbers and polynomials. Furthermore, we
  construct Euler $q$-Zeta function.
\endabstract
%\thanks  2000 AMS Subject Classification:  11B68, 11S40.
%\newline keywords and phrases :Sums of powers, Bernoulli numbers, $q$-Bernoulli numbers, zeta
%function, Dirichlet series
%\newline $\ast$ This work was supported by grant No. (R08-2003-10439-0) from the
%Basic Research Program of the Korea Science \& Engineering
%Foundation.
%\endthanks \NoRunningHeads
\endtopmatter

\document

\head{ 1. Introduction }
\endhead

For any complex number $z$, it is well known that the familiar
Euler polynomials $E_n (z)$ are defined by means of the following
generating function, cf.[3, 5, 6, 9, 13]:

$$\eqalignno{&
F(z,t)= \frac{2}{e^t  +1} e^{zt}  = \sum_{n=0}^{\infty} E_n (z)
\frac{t^n}{ n!},\ \  (|t|<\pi ).&(1)}
$$

We note that, by substituting $z=0$ into (1), $E_n (0) = E_n$ is
the familiar $n$-th Euler number defined by
$$
G(t) =F(0,t) =\frac{2}{e^t  +1} =
 \sum_{n=0}^{\infty} E_n
\frac{t^n}{ n!},\ \  (|t|<\pi ), \text{ cf.[4, 5].}
$$

By the meaning of the generalization of $E_n$, Frobenius-Euler
numbers and polynomials are also defined by

 $$\frac{1-u}{e^t  -u}
= \sum_{n=0}^{\infty} H_n (u) \frac{t^n}{n!},\text{ and
}\frac{1-u}{e^t -u} e^{xt} =\sum_{n=0}^{\infty} H_n (u,x )
\frac{t^n}{n!}\ \ (u\in \Bbb C \text{ with } |u|>1), \text{
cf.[16].}$$

Over five decades ago, Calitz [2, 3] defined $q$-extension of
Frobenius-Euler  numbers and polynomials and proved properties
analogous to those satisfied $H_n (u)$ and $H_n (u,x)$.  Recently,
Satoh [14, 15] used these properties, especially the so-called
distribution relation for the $q$-Frobenius-Euler polynomials, in
order to construct the corresponding $q$-extension of the $p$-adic
measure and to define a $q$-extension of $p$-adic $l$-function
$l_{p,q} (s,u)$.

Let $p$ be a fixed odd prime in this paper. Throughout this paper,
the symbols $\Bbb Z$, $\Bbb Z_p$, $\Bbb Q_p$, $\Bbb C$ and $\Bbb
C_p$, denote the ring of rational integers, the ring of $p$-adic
integers, the field of $p$-adic numbers, the complex number field,
and the completion of the algebraic closure of $\Bbb Q_p$,
respectively. Let $\nu_p (p)$ be the normalized exponential
valuation of $\Bbb C_p$ with $|p|_p = p^{-\nu_p (p)} =p^{-1}$.
When one speaks of $q$-extension, $q$ can be regarded as an
indeterminate, a complex number $q\in \Bbb C$, or a $p$-adic
number $q\in \Bbb C_p$; it is always clear  from the context. If
$q\in \Bbb C$, then one usually assumes that $|q|<1$. If $q\in
\Bbb C_p$, then one usually assumes that $|q-1|_p
<p^{-\frac{1}{p-1}}$, and hence $q^x =\exp (x \log q )$ for $x\in
\Bbb Z_p $. In this paper, we use the below notation
$$[x]_q =\frac{1-q^x}{ 1-q}, \ \ \  (a:q)_n =
(1-a)(1-aq)\cdots(1-aq^{n-1}), \text{ cf.[6, 7, 8, 9, 10, 11, 14
].}$$ Note that $\lim_{q\to 1} [x]_q =x$ for any $x$ with $|x|_p
\leq 1$ in  the $p$-adic case. For a fixed positive integer $d$
with $(p,d)=1$, set
$$\split &
X =X_d =\varprojlim_N \Bbb Z/dp^N ,\cr &  X_1 = \Bbb Z_p , X^\ast
=\bigcup_{\Sb  0<a <dp\\(a,p)=1
\endSb} a +d p\Bbb Z_p ,\\
& a +d p^N \Bbb Z_p = \{ x\in X | x \equiv a
\pmod{p^N}\},\endsplit
$$
where $a\in\Bbb Z$ satisfies the condition $0\leq a < d p^N$,
(cf.[10,11]). We say that $f$ is a uniformly differentiable
function at a point $a\in\Bbb Z_p$, and write $f\in UD(\Bbb Z_p
)$, if the difference quotients $F_f (x,y) = \frac{f(x)
-f(y)}{x-y}$ have a limit $f^\prime (a)$ as $(x,y)\to (a,a)$,
cf.[11]. For $f\in UD (\Bbb Z_p )$, let us begin with the
expression
$$
\frac{1}{[p^N ]_q } \sum_{0\leq j< p^N } q^j f(j) =\sum_{0\leq j<
p^N } f(j) \mu_q (j+ p^N \Bbb Z_p ), \text{ cf.[7, 8, 9, 11],}
$$
which represents  a $q$-analogue of Riemann sums for $f$. The
integral of $f$ on $\Bbb Z_p$ is defined as the limit of those
sums(as $n\to \infty$) if this limit exists. The $q$-Volkenborn
integral of a function $f \in UD(\Bbb Z_p )$ is defined by
$$I_q (f) = \int_X f(x) d\mu_q (x) =\int_{X_d} f(x) d\mu_q (x)=
\lim_{N\to\infty} \frac{1}{[dp^N]_q} \sum_{x=0}^{dp^N -1}f(x) q^x
.$$ Recently, we considered another construction of a $q$-Eulerian
numbers, which are different than Carlitz's $q$-Eulerian numbers
as follows, cf.[6, 12, 13]:
$$
F_q (x,t)=[2]_q \sum_{n=0}^\infty (-1)^n q^n e^{[n+x]_q t}
=\sum_{n=0}^\infty E_{n,q}(x) \frac{t^n}{n!}.
$$
Thus  we have
$$\split &
E_{n,q} =E_{n,q}(0) =\frac{[2]_q}{(1-q)_n} \sum_{l=0}^n
\binom{n}{l} \frac{(-1)^l }{1+q^{l+1}},\cr & E_{n,q}(x)
=\frac{[2]_q}{(1-q)_n} \sum_{l=0}^n \binom{n}{l} \frac{(-1)^l
}{1+q^{l+1}} q^{lx},
\endsplit
$$
where $\binom{n}{l}$ is a binomial coefficient, cf.[13].

Note that $\lim_{q\to 1} E_{n,q} = E_n$ and $\lim_{q\to 1}
E_{n,q}(x) = E_n (x)$. In [12], we also proved that $q$-Eulerian
polynomial $E_{n,q} (x)$ can be  represented by $q$-Volkenborn
integral as follows:
$$
\int_{X_d} [x+x_1 ]_q^k d\mu_{-q} (x_1 ) =\int_{\Bbb Z_p} [x+x_1
]_q^k d\mu_{-q} (x ) =E_{k,q} (x) , \ \ \text{ for } k,d\in \Bbb
N,
$$
where $\mu_{-q} (x+ p^N \Bbb Z_p ) =\frac{q^x [2]_q }{1+q^{p^N}}
(-1)^x .$

The purpose of this paper is to present a systemic study of some
families of multiple $q$-Euler numbers and polynomials. In
particular, by using the $q$-Volkenborn integration on $\Bbb Z_p$,
we construct $p$-adic $q$-Euler numbers and polynomials of higher
order. We also define new generating function of these $q$-Euler
numbers and polynomials of higher order. Furthermore, we construct
Euler $q$-$\zeta$-function. From section 2 to section 5, we assume
that $q\in\Bbb C_p$ with $|1-q|_p < p^{-\frac{1}{p-1}}$.

\head{ 2. $q$-Euler numbers and polynomials associated with an
invariant $p$-adic $q$-integrals on $\Bbb Z_p$ }
\endhead

Let $h\in \Bbb Z$, $k\in\Bbb N=\{1,2,3,\cdots\}$, and let us
consider the extended higher-order $q$-Euler numbers as follows:
$$
E_{m,q}^{(h,k)} =\underbrace{\int_{\Bbb Z_p}\cdots\int_{\Bbb
Z_p}}_{k \ \text{times} } [x_1 +x_2 +\cdots +x_k ]_q^m q^{x_1
(h-1) +\cdots+ x_k (h-k)} d\mu_{-q} (x_1 )\cdots d\mu_{-q} (x_k )
.$$

Then we have
$$
E_{m,q}^{(h,k)} = \frac{[2]_q^k}{ (1-q)^m }\sum_{l=0}^m
\binom{m}{l} \frac{(-1)^l}{(-q^{h+l }: q^{-1})_k} .
$$
From the definition of $E_{m,q}^{(h , k)}$, we can easily derive
the below:
$$E_{m,q}^{(h,k)} = E_{m,q}^{(h-1 ,k) } + (q-1)
E_{m+1,q}^{(h-1,k)}, \ \ (m \geq 0).$$ It is easy to show that
$$\eqalignno{&
\underbrace{\int_{\Bbb Z_p}\cdots\int_{\Bbb Z_p}}_{k+1 \
\text{times} } q^{\sum_{j=1}^{k+1} (m-j)x_j} d\mu_{-q} (x_1 )
\cdots d\mu_{-q} (x_{k+1} )\cr & = \sum_{j=1}^m \binom{m}{j}
(q-1)^j \underbrace{\int_{\Bbb Z_p}\cdots\int_{\Bbb Z_p}}_{k+1 \
\text{times} } [ \sum_{l=1}^{k+1} x_l ]_q^j \  q^{- \sum_{j=1}^k
jx_j} d\mu_{-q} (x_1 )\cdots d\mu_{-q} (x_{k+1} ) ,&(2)}
$$
and we also get
$$\eqalignno{&
\underbrace{\int_{\Bbb Z_p}\cdots\int_{\Bbb Z_p}}_{k+1 \
\text{times} } q^{\sum_{j=1}^{k+1} (m-j)x_j} d\mu_{-q} (x_1 )
\cdots d\mu_{-q} (x_{k+1} )= \frac{[2]_q^{k+1}}{ (-q^m :
q^{-1})_{k+1}} .&(3)}$$ From (2) and (3), we can derive the below
proposition.
\proclaim{Proposition 1} For $m,k \in\Bbb N$, we have
$$
\split & \sum_{j=0}^m \binom{m}{j} (q-1)^j
E_{j,q}^{(0,k+1)}=\frac{[2]_q^{k+1} }{ (-q^m : q^{-1})_{k+1}} ,\cr
 &
E_{m,q}^{(h,k)} = \frac{[2]_q^k}{ (1-q)^m} \sum_{l=0}^m
\binom{m}{l} \frac{(-1)^l }{ (-q^{h+l} : q^{-1} )_k } .
\endsplit
$$
\endproclaim
\remark{Remark} Note that $E_{n,q}^{(1,1)} =E_{n,q}$, where
$E_{n,q}$ are the $q$-Euler numbers(see [13]).
\endremark
 From the definition
of $E_{n,q}^{(h,k)}$, we can derive
$$
\sum_{j=0}^i \binom{i}{j} (q-1)^j E_{m-i+j,q}^{(h-1 ,k)}
=\sum_{j=0}^{i-1} \binom{i-1}{j} (q-1)^j E_{m+j-i,q}^{(h,k)}
$$
for $m\geq i$. By simple calculation, we easily see that
$$
\sum_{j=0}^m \binom{m}{j} (q-1)^j E_{j,q}^{(h ,1)} =\int_{\Bbb
Z_p} q^{mx} q^{(h-1)x} d\mu_{-q}(x) = \frac{[2]_q}{[2]_{q^{m+h}}}.
$$
Furthermore, we can give the following relation for the $q$-Euler
numbers, $E_{m,q}^{(0,h)}$,:
$$\eqalignno{&
\sum_{j=0}^m \binom{m}{j} (q-1)^j E_{j,q}^{(0 ,k)} = \frac{[2]_q^k
}{(-q^m : q^{-1})_k}.&(4)}
$$
\head{ 3. Polynomials $E_{n,q}^{(0,k)} (x)$ }
\endhead
We now define
 the polynomials $E_{n,q}^{(0,k)} (x)$ (in $q^x$) by
 $$
 E_{n,q}^{(0,k)} (x) = \underbrace{\int_{\Bbb Z_p}\cdots\int_{\Bbb
Z_p}}_{k \ \text{times} }  [x_1 +x_2 +\cdots +x_k ]_q^m
q^{\sum_{j=1}^k j x_j } d\mu_{-q} (x_1 )\cdots d\mu_{-q} (x_k ) .
 $$
Thus, we have
$$\eqalignno{&
(q-1)^m E_{m,q}^{(0,k)} (x) =[2]_q^k
 \sum_{j=0}^m \binom{m}{j} q^{jx} (-1)^{m-j} \frac{1}{(-q^j : q^{-1})_k }
. &(5)}
$$
It is not difficult to show that
$$\underbrace{\int_{\Bbb Z_p}\cdots\int_{\Bbb Z_p}}_{k \
\text{times} } q^{\sum_{j=1}^m (m-j)x_j +mx} d\mu_{-q} (x_1
)\cdots d\mu_{-q} (x_k ) = q^{mx} \frac{ [2]_q^k}{(-q^m :
q^{-1})_k},$$ and
$$
\underbrace{\int_{\Bbb Z_p}\cdots\int_{\Bbb Z_p}}_{k \
\text{times} } q^{\sum_{j=1}^m (m-j)x_j +mx} d\mu_{-q} (x_1
)\cdots d\mu_{-q} (x_k ) =\sum_{j=0}^m \binom{m}{j} (q-1)^j
 E_{j,q}^{(0,k)} (x) .$$
 Therefore we obtain the following.

\proclaim{Lemma 2} For $m,k \in\Bbb N$, we have
$$\eqalignno{&
 \sum_{j=0}^m \binom{m}{j} (q-1)^j E_{j,q}^{(0,k)} (x)
=\frac{q^{mx}[2]_q^k }{ (-q^m : q^{-1})_{k}} ,\cr
 &
E_{m,q}^{(0,k)}(x) = \frac{[2]_q^k}{ (1-q)^m} \sum_{j=0}^m
\binom{m}{j} q^{jx} (-1)^j \frac{1 }{ (-q^{j} : q^{-1} )_k } .
&(6)}
$$
\endproclaim

Let $l\in\Bbb N$ with $l\equiv 1 \pmod2$. Then we get easily
$$\split &
\underbrace{\int_{\Bbb Z_p}\cdots\int_{\Bbb Z_p}}_{k \
\text{times} } \left[ x+\sum_{j=1}^k x_j \right]_q^m
q^{-\sum_{j=1}^k j x_j } d\mu_{-q} (x_1 )\cdots d\mu_{-q} (x_k
)\cr &= \frac{[l]_q^m}{[l]_{-q}^k} \sum_{i_1 , \cdots, i_k
=0}^{l-1} q^{-\sum_{j=2}^k (j-1)i_j } \cr &\ \ \
\cdot(-1)^{\sum_{j=1}^k i_j} \underbrace{\int_{\Bbb
Z_p}\cdots\int_{\Bbb Z_p}}_{k \ \text{times} } \left[ \frac{x+
\sum_{j=1}^k i_j }{l}+ \sum_{j=1}^k x_j \right]_{q^l}^m  q^{-l
\sum_{j=1}^k j x_j} d\mu_{-q^l} (x_1 )\cdots d\mu_{-q^l} (x_k ).
\endsplit$$
From this, we can derive the following ``multiplication formula":
\proclaim{Theorem 3} Let  $l$ be an odd positive integer. Then
$$\eqalignno{ &  E_{m,q}^{(0,k)} (x) =
\frac{[l]_q^m}{[l]_{-q}^k}  \sum_{i_1 , \cdots, i_k =0}^{l-1}
q^{-\sum_{j=2}^k (j-1)i_j } (-1)^{\sum_{l=1}^k i_l} E_{m,q^l
}^{(0,k)} (\frac{x+i_1 +\cdots+i_k}{l}). &(7)}
$$
Moreover,
$$\eqalignno{ &  E_{m,q}^{(0,k)} (lx) =
\frac{[l]_q^m}{[l]_{-q}^k}  \sum_{i_1 , \cdots, i_k =0}^{l-1}
q^{-\sum_{j=2}^k (j-1)i_j } (-1)^{\sum_{l=1}^k i_l} E_{m,q^l
}^{(0,k)} (x+ \frac{i_1 +\cdots+i_k}{l}). &(8)}
$$\endproclaim
From (4) and (5), we can also derive the below expression for
$E_{n,q}^{(0,k)} (x)$:
$$\eqalignno{ &  E_{m,q}^{(0,k)} (x) =
  \sum_{i=0}^{m} \binom{m}{i} E_{i,q}^{(0,k)} [x]_q^{m-i} q^{ix}, &(9)}$$
whence also
$$\eqalignno{ &  E_{m,q}^{(0,k)} (x+y) =
  \sum_{j=0}^{m} \binom{m}{j}  [y]_q^{m-i} q^{jy} E_{j,q}^{(0,k)}(x) . &(10)}$$

\head{ 4.  Polynomials  $E_{m,q}^{(h,1)}(x)$ }
\endhead

Let us define
$$\eqalignno{ &
E_{m,q}^{(h,1)}(x) =\int_{\Bbb Z_p} [x+x_1 ]_q^m \ q^{x_1
(h-1)}d\mu_{-q} (x_1 ). &(11)}
$$
Then we have
$$E_{m,q}^{(h,1)}(x) = \frac{[2]_q }{(1-q)^m} \sum_{l=0}^m
\binom{m}{l} (-1)^l q^{lx} \frac{1}{(1+q^{l+h})}.
$$

By simple calculation of $q$-Volkenvorn integral, we note that
$$\split
q^x \int_{\Bbb Z_p} [x+x_1 ]_q^m \ q^{x_1 (h-1)}d\mu_{-q} (x_1 )
&= (q-1) \int_{\Bbb Z_p} [x+x_1 ]_q^{m+1} \ q^{x_1 (h-2)}d\mu_{-q}
(x_1 )\cr & +
 \int_{\Bbb Z_p} [x+x_1 ]_q^{m} \ q^{x_1 (h-2)}d\mu_{-q}
(x_1 ).\endsplit
$$
Thus, we have
$$\eqalignno{ &
q^x E_{m,q}^{(h,1)}(x) =(q-1)E_{m+1,q}^{(h-1,1)}(x) +
E_{m,q}^{(h-1,1)}(x). &(12)}
$$
It is easy to show that
$$\int_{\Bbb Z_p} [x+x_1 ]_q^m \ q^{ (h-1)x_1} d\mu_{-q} (x_1 )
=\sum_{j=0}^m \binom{m}{j} [x]_q^{m-j} q^{jx} \int_{\Bbb Z_p} [x_1
]_q^j q^{(h-1)x_1}d\mu_{-q} (x_1 ).
$$
This is equivalent to
$$\split
E_{m,q}^{(h,1)}(x) &=\sum_{j=0}^m \binom{m}{j} [x]_q^{m-j} q^{jx}
E_{j,q}^{(h,1)} \cr &=\left( q^x E_q^{(h,1)} +[x]_q \right)^m , \
\ \text{ for } m\geq 1,\endsplit$$ where we use the technique
method notation by replacing $(E_q^{(h,1)})^n$ by
$E_{n,q}^{(h,1)}$, symbolically. From (11), we can derive
$$\eqalignno{ &
q^h E_{m,q}^{(h,1)}(x+1) +E_{m,q}^{(h,1)}(x)=[2]_q [x]_q^m .
&(13)}
$$

For $x=0$ in (13), this gives
$$\eqalignno{ &
q^h \left( q E_{m,q}^{(h,1)}+1\right)^m
+E_{m,q}^{(h,1)}=\delta_{0,k}, &(14)}
$$
where $\delta_{0,k}$ is Kronecker symbol. By the simple
calculation of $q$-Volkenborn integration, we easily see that
$$ \int_{\Bbb Z_p}  q^{x_1 (h-1)}d\mu_{-q}
(x_1 ) = \frac{[2]_q}{[2]_{q^h}}.$$ Thus, we have
$E_{0,q}^{(h,1)}=\frac{[2]_q}{[2]_{q^h}}$. From the definition of
$q$-Euler polynomials, we can derive
$$\int_{\Bbb Z_p} [1-x+x_1 ]_{q-1}^m \ q^{-x_1 (h-1)} d\mu_{-q} (x_1
)= q^{m+h-1} (-1)^m E_{m,q}^{(h,1)} (x).$$ Therefore we obtain the
below ``complementary formula": \proclaim{Theorem 4} For $m\in\Bbb
N$, $n\in \Bbb Z$, we have
$$\eqalignno{ &  E_{m,q^{-1} }^{(h,1)} (1-x) =(-1)^m q^{m+h-1}E_{m,q}^{(h,1)} (x). &(15)}
$$
In particular, for $x=1$, we see that
$$\eqalignno{   E_{m,q-1}^{(h,1)} (0) &=(-1)^m q^{m+h-1}E_{m,q}^{(h,1)} (1)\cr
&= (-1)^{m-1} q^{m-1} E_{m,q}^{(h,1)}, \ \ \text{ for } m\geq 1  .
&(16)}
$$
\endproclaim

For $l\in \Bbb N$ with $l\equiv 1\pmod 2$, we have
$$\split &
\int_{\Bbb Z_p}q^{(h-1)x_1 }  [x+x_1 ]_{q}^m \ q^{x_1 (h-1)}
d\mu_{-q} (x_1 )\cr &= \frac{[l]_q^m}{[l]_{-q} } \sum_{i=0}^{l-1}
q^{hi} (-1)^i \int_{\Bbb Z_p} \left[ \frac{x+i}{l}+x_1
\right]_{q^l}^m \ q^{x_1 (h-1)l} d\mu_{-q^l} (x_1 ).\endsplit$$
Thus, we can also obtain the following:
 \proclaim{Theorem
5}(Multiplication formula) For $l\in\Bbb N$ with $l\equiv 1
\pmod2$, we have
$$\frac{[2]_q}{[2]_{q^l}} [l]_q^m \sum_{i=0}^{l-1}
q^{hi} (-1)^i  E_{m,q^l}^{(h,1)} (\frac{x+i }{l}) =
E_{m,q}^{(h,1)} (x).$$
Furthermore,
$$\frac{[2]_q}{[2]_{q^l}} [l]_q^m \sum_{i=0}^{l-1}
q^{hi} (-1)^i  E_{m,q^l}^{(h,1)} (x+ \frac{i }{l})
=E_{m,q}^{(h,1)} (lx).
$$
\endproclaim

\head{ 5. Polynomials $E_{m,q}^{(h,k)} (x)$  and $h=k$}
\endhead

It is now easy to combine the above results and define the new
polynomials as follows:

$$
E_{m,q}^{(h,k)} (x) =\underbrace{\int_{\Bbb Z_p}\cdots\int_{\Bbb
Z_p}}_{k \ \text{times} } [x+x_1 +\cdots+x_k]_q^m
 q^{(h-1)x_1 +\cdots+ (h-k)x_k } d\mu_{-q}
(x_1 )\cdots d\mu_{-q} (x_k ).
$$
Thus, we note that
$$\eqalignno{ &
(q-1)^m  E_{m,q}^{(h,k)} (x) =\sum_{j=0}^m \binom{m}{j} (-1)^{m-j}
q^{xj} \frac{[2]_q^k}{(-q^{j+h}:q^{-1})_k}  . &(17)}$$ We may now
mention the following formulas which are easy to prove.
$$\eqalignno{ &
q^h  E_{m,q}^{(h,k)} (x+1)+E_{m,q}^{(h,k)} (x) =[2]_q
E_{m,q}^{(h-1,k-1)} (x), &(18)}
$$
 and
$$\eqalignno{ &
q^x  E_{m,q}^{(h+1,k)} (x) =(q-1) E_{m+1,q}^{(h,k)}
(x)+E_{m,q}^{(h,k)} (x). &(19)}$$ Let $l\in \Bbb N$ with $l\equiv
1\pmod 2$. Then  we note that
$$\split &
\underbrace{\int_{\Bbb Z_p}\cdots\int_{\Bbb Z_p}}_{k \
\text{times} } [x+\sum_{j=1}^k x_j ]_q^m
 q^{\sum_{j=1}^k(h-j)x_j } d\mu_{-q}
(x_1 )\cdots d\mu_{-q} (x_k )\cr &= \frac{[l]_q^m}{[l]_{-q}^k}
\sum_{i_1,\cdots, i_k =0}^{l-1} q^{ h\sum_{j=1}^k i_j
-\sum_{j=2}^k (j-1)i_j  } (-1)^{\sum_{j=1}^k i_j}\cr &\cdot
\underbrace{\int_{\Bbb Z_p}\cdots\int_{\Bbb Z_p}}_{k \
\text{times} } \left[ \frac{x+\sum_{j=1}^k i_j}{l} +\sum_{j=1}^k
x_j \right]_{q^l}^m
 (q^l)^{\sum_{j=1}^k (h-j)x_j } d\mu_{-q^l}
(x_1 )\cdots d\mu_{-q^l } (x_k ).\endsplit$$ Therefore we obtain
the following:

\proclaim{Theorem 6}( Distribution for q-Euler polynomials)
 For $l\in\Bbb N$ with $l\equiv
1 \pmod2$. Then we have
$$E_{m,q}^{(h,k)} (lx) = \frac{[l]_q^m}{[l]_{-q}^k}
\sum_{i_1,\cdot,i_k=0}^{l-1} q^{h\sum_{j=1}^k i_j -\sum_{j=2}^k
(j-1)i_j} (-1)^{\sum_{j=1}^{k} i_j} E_{m,q^l}^{(h,k)} \left( x+
\frac{i_1 +\cdots+i_k }{l}\right).\tag20$$
\endproclaim

It is interesting to consider the case $h=k$, which leads to the
desired extension of the $q$-Euler numbers of higher order,
cf.[1]. We shall denote the polynomials in this special case by
$E_{m,q}^{(k)}(x):=E_{m,q}^{(k,k)}(x)$. Then we have

$$\eqalignno{ &
(q-1)^m E_{m,q}^{(k)} (x) = \sum_{j=0}^m \binom{m}{j} (-1)^{m-j}
q^{jx}
 \frac{[2]_q^k}{(-q^{j+k}: q^{-1})_k}, &(21)}
$$
and
$$\eqalignno{ &
 E_{m,q^{-1}}^{(k)} (k-x) =(-1)^m q^{m+\binom{k}{2}}E_{m,q}^{(k)} (x). &(22)}$$
For $x=k$ in (22), we see that
$$\eqalignno{ &
 E_{m,q^{-1}}^{(k)} (0) =(-1)^m q^{m+\binom{k}{2}}E_{m,q}^{(k)} (k). &(23)}$$
From (18), we can derive the below formula:
$$\eqalignno{ &
q^k  E_{m,q}^{(k)} (x+1) + E_{m,q}^{(k)} (x)=[2]_q E_{m,q}^{(k-1)}
(x). &(24)}
$$
Putting $x=0$ in (17), we obtain
$$\eqalignno{ &
(q-1)^m  E_{m,q}^{(k)} =\sum_{i=0}^m \binom{m}{i}
(-1)^{m-i}\frac{[2]_q^k}{(-q^{i+k}:q^{-1})_k}. &(25)}
$$
Note that
$$\split &
\sum_{i=0}^m  \binom{m}{i}(q-1)^i \underbrace{\int_{\Bbb
Z_p}\cdots\int_{\Bbb Z_p}}_{k \ \text{times} } [x_1 +\cdots+x_k
]_q^i
 q^{\sum_{j=1}^{k-1} (k-j)x_j } d\mu_{-q}
(x_1 )\cdots d\mu_{-q} (x_k )\cr &= \frac{[2]_q^k}{(-q^{m+k}:
q^{-1})_k}.\endsplit
$$

From this, we can easily derive
$$\eqalignno{ &
\sum_{i=0}^m \binom{m}{i} (q-1)^i E_{i,q}^{(k)} =
\frac{[2]_q^k}{(-q^{m+k}:q^{-1})_k} &(26)}
$$
and so it follows
$$\eqalignno{ &E_{m,q}^{(k)} (x) = (q^x E_{q}^{(k)} +[x]_q)^m , \
\ m\geq 1, &(27)}
$$
where we use the technique method notation by replacing
$(E_q^{(k)})^n$ by $E_{n,q}^{(k)}$, symbolically. In particular,
from (24), we have
$$\eqalignno{ & q^k (q E_{q}^{(k)} +1)^m +E_{m,q}^{(k)} =[2]_q E_{m,q}^{(k-1)}.  &(28)}
$$

It is easy to see that
$$
\underbrace{\int_{\Bbb Z_p}\cdots\int_{\Bbb Z_p}}_{k \
\text{times} }
 q^{(k-1)x_1 +\cdots+x_{k-1}} d\mu_{-q}
(x_1 )\cdots d\mu_{-q} (x_k ) =\frac{[2]_q^k}{(-q^k :q^{-1})_k}.
$$

Thus, we note that $E_{0,q}^{(k)} =\frac{[2]_q^k}{(-q^k
:q^{-1})_k}$.

\head{ 6. Generating function for $q$-Euler polynomials }
\endhead

An obvious generating function for $q$-Euler polynomials is
obtained, from  (17), by
$$\eqalignno{ & [2]_q^k e^{\frac{t}{1-q}} \sum_{j=0}^\infty
\frac{(-1)^j}{(-q^{j+h}: q^{-1})_k} q^{jx} \left( \frac{1}{1-q}
\right)^j \frac{t^j}{j!}\cr &= \sum_{n=0}^\infty E_{n,q}^{(h,k)}
\frac{t^n}{n!}.&(29)}
$$

From (17), we can also derive the below formula:

$$\eqalignno{ & q^{h-k}
E_{m,q}^{(h,k+1)} (x+1) =[2]_q E_{m,q}^{(h,k)}(x) -
E_{m,q}^{(h,k+1)}(x).&(30)}
$$

Again from (21) and (25), we get easily
$$\split &
\underbrace{\int_{\Bbb Z_p}\cdots\int_{\Bbb Z_p}}_{k \
\text{times} } [x+\sum_{j=1}^k x_j]_q^m
 q^{\sum_{j=1}^{k-1} (k-j)x_j} d\mu_{-q}
(x_1 )\cdots d\mu_{-q} (x_k )\cr & = \sum_{j=0}^m
\binom{m}{j}q^{xj} \underbrace{\int_{\Bbb Z_p}\cdots\int_{\Bbb
Z_p}}_{k \ \text{times} } [x_k]_q^j [x+\sum_{j=1}^{k-1}
x_j]_q^{n-j}
 q^{\sum_{l=1}^{k-1} (k+j-l)x_l} d\mu_{-q}
(x_1 )\cdots d\mu_{-q} (x_k ).\endsplit$$ Thus, we note that
$$\eqalignno{ &
E_{m,q}^{(k)} (x) = \sum_{j=0}^m \binom{m}{j} q^{xj} E_{j,q}^{(1)}
 E_{m-j,q}^{(k+j,k-1)}(x).&(31)}
$$
Take $x=0$ in (31), we have
$$\eqalignno{ &
E_{m,q}^{(k)}  = \sum_{i=0}^m \binom{m}{i}  E_{j,q}^{(1)}
 E_{m-j,q}^{(k+j,k-1)}.&(32)}$$
So, for $k=2$,
$$E_{m,q}^{(2)}  = \sum_{i=0}^m \binom{m}{i}  E_{j,q}
 E_{m-j,q}^{(j+2,1)}.$$
It is not difficult to show that
$$\int_{\Bbb Z_p} [x]_q^m q^{hx} d\mu_{-q} (x) =\sum_{j=0}^h
\binom{h}{j} (q-1)^j \int_{\Bbb Z_p} [x]_q^{m+j}  d\mu_{-q} (x),
\text{ for $h\in \Bbb N$ }.$$  From this, we can derive the below:
$$\eqalignno{ &
E_{m,q}^{(h+1,1)}  = \sum_{j=0}^h \binom{h}{j} (q-1)^j E_{m+j,q},
\ \ h\in \Bbb N .&(33)}$$  By (32) and (33), we easily see that
$$\eqalignno{ &
E_{m,q}^{(2)}  = \sum_{j=0}^m \binom{m}{j}
E_{j,q}\sum_{i=0}^{j+1}\binom{j+1}{i} (q-1)^i E_{m-j+i,q} .&(34)}
$$

By (34), for $q=1$, we note that
$$E_m^{(2)} = \sum_{j=0}^m \binom{m}{j} E_j E_{m-j},
\text{ where $\left(\frac{2}{e^t +1}\right)^k = \sum_{n=0}^\infty
E_n^{(k)} \frac{t^n}{n!}$.}$$  It is easy to show that
$$
[x+x_1 +\cdots+x_k ]_q^m =\sum_{j=0}^m \binom{m}{j} [x_1
+x]_q^{m-j} q^{j(x_1 +x)} [x_2 +\cdots+ x_k ]_q^j .
$$
By using this, we get easily
$$\split &
\underbrace{\int_{\Bbb Z_p}\cdots\int_{\Bbb Z_p}}_{k \
\text{times} } [x+\sum_{j=1}^k x_j]_q^m
 q^{\sum_{j=1}^{k-1} (k-j)x_j} d\mu_{-q}
(x_1 )\cdots d\mu_{-q} (x_k ) \cr &= \sum_{j=0}^m \binom{m}{j}
q^{jx}
 \int_{\Bbb Z_p} [x+x_1]_q^{m-j} q^{(k+j-1)x_1} d\mu_{-q}
(x_1 )\cr &\quad\cdot \underbrace{\int_{\Bbb Z_p}\cdots\int_{\Bbb
Z_p}}_{k-1 \ \text{times} } [x_2 +\cdots+ x_k]_q^j
 q^{\sum_{j=2}^{k-1} (k-j)x_j} d\mu_{-q}
(x_2 )\cdots d\mu_{-q} (x_k ).\endsplit
$$

Therefore we obtain the following:

\proclaim{Theorem 7} For $m,k\in\Bbb N$, we have
$$\eqalignno{ &
E_{m,q}^{(k)} (x) =\sum_{j=0}^m \binom{m}{j} q^{jx}
E_{m-j,q}^{(k+j,1)}(x)E_{j,q}^{(k-1)}.&(35)}
$$

Indeed for $x=0$,
$$\eqalignno{ &
E_{m,q}^{(k)}  =\sum_{j=0}^m \binom{m}{j} E_{m-j,q}^{(k+j,1)}
E_{j,q}^{(k-1)}\cr &= \sum_{j=0}^m \binom{m}{j} E_{j,q}^{(k-1)}
\sum_{j=0}^{k+j} (q-1)^i \binom{k+j-1}{i} E_{m-j+i,q}^{(1)}
.&(36)}
$$
As for  $q=1$, we get the below formula
$$E_{m}^{(k)}  =\sum_{j=0}^m \binom{m}{j} E_{j}^{(k-1)}
E_{m-j}^{(1)}.
$$
\endproclaim

\head{ 7.  $q$-Euler zeta function in $\Bbb C$ }
\endhead

In this section, we assume that $q\in\Bbb C$ with $|q|<1$. From
section 4, we note that

$$\eqalignno{
E_{m,q}^{(h,1)} (x) & = \frac{[2]_q}{(1-q)^m} \sum_{l=0}^m
\binom{m}{l}   q^{lx} (-1)^l \frac{1}{1+q^{l+h}}\cr
 &= [2]_q  \sum_{n=0}^\infty (-1)^n q^{nh} [n+x]_q^n .&(37)}
$$
Thus, we can define $q$-Euler zeta function:

\definition{Defintion 8} For $s,q\in\Bbb C$ with $|q|<1$, define
$$
\zeta_{E,q}^h (s,x) =[2]_q \sum_{n=0}^\infty \frac{(-1)^n q^{nh}
}{[n+x]_q^s},
$$
where $x\in \Bbb R$ with $0< x\leq 1$.
\enddefinition

Note that $\zeta_{E,q}^h (-m,x) =E_{m,q}^{(h,1)} (x)$, for $m\in
\Bbb N$. Let
$$
F_q (t,x ) =\sum_{n=0}^\infty E_{n,q}^{(h,1)} (x) \frac{t^n }{n!}.
$$
Then we have
$$\split
F_q (t,x )& = [2]_q e^{\frac{t}{1-q}} \sum_{n=0}^\infty  (-1)^n
q^{hn} e^{-\frac{q^{n+x}}{1-q}t }\cr
 &= [2]_q \sum_{n=0}^\infty  (-1)^n
q^{hn} e^{[n+x]_q t}, \ \text{ for } h\in\Bbb Z .
 \endsplit$$

Therefore we obtain the following

\proclaim{Lemma 9} For $h\in\Bbb Z$, we have
$$\eqalignno{
F_q (t,x )& = [2]_q \sum_{n=0}^\infty  (-1)^n q^{hn}  e^{[n+x]_q
t}\cr   & =\sum_{n=0}^\infty  E_{n,q}^{(h,1)} (x) \frac{t^n
}{n!}.&(38)}
$$

\endproclaim

Let $\Gamma (s)$ be the gamma function. Then we easily see that

$$\eqalignno{&
\frac{1}{\Gamma (s)} \int_0^\infty t^{s-1} F_q (-t, x) dt =
\zeta_{E,q}^h (s,x ), \ \ \text{ for }  s\in \Bbb C.&(39)}
$$

From (38) and (39), we can also  derive the below Eq.(40):

$$
\zeta_{E,q}^h (-n, x) =E_{n,q}^{(h,1)} (x), \ \ \text{ for }
n\in\Bbb N. \tag40$$

    \enddocument